\documentclass{article}
\usepackage[T1]{fontenc}
\usepackage[utf8]{inputenc}
\usepackage{amsthm}
\usepackage{amssymb}
\usepackage{amsmath}
\usepackage{graphicx}
\usepackage{float}
\usepackage{caption}
\usepackage{url}
\usepackage{hyperref}
\usepackage{enumitem}
\usepackage{mathtools}
\theoremstyle{definition}

\theoremstyle{remark}

\author{M.~Hellus
\and A.~Rechenauer
\and R.~Waldi}
\title{A lower bound for the Wilf density, deduced from a result of Zhai}
\begin{document}
\maketitle
\begin{abstract}
Let $S\neq\mathbb N$ be a numerical semigroup with Frobenius number $f$, genus $g$ and embedding dimension $e$. In 1978 Wilf asked the question, whether $\frac{f+1-g}{f+1}\geq\frac1e$. As is well known, this holds in the cases $e=2$ and $e=3$.

From Zhai's results in \cite{Zhai} we derive
\[\frac{f+1-g}{f+1}\geq\frac2{e^2-e+2}\text{ for }e\geq4\,.\]
\end{abstract}
Let $\mathbb N$ be the set of nonnegative integers. A \emph{numerical semigroup} is an additively closed subset $S$ of $\mathbb N$ with $0\in S$ and only finitely many positive integers outside $S$, the so-called \emph{gaps} of $S$. We suppose $S$ to be a numerical semigroup different from $\mathbb N$. The number of gaps is called the \emph{genus} $g$, the largest gap is the \emph{Frobenius} number $f$ of $S$. The set $E=S^*\setminus(S^*+S^*)$, where $S^*= S\setminus\{0\}$, is the unique minimal system of generators of $S$. The latter are called the \emph{atoms} of $S$; their number $e$ is the \emph{embedding dimension} of $S$, the smallest atom is the \emph{multiplicity} $m$ of $S$.

In 1978, Wilf \cite{Wilf} posed the following question: Does one always have
\[\tag{1}d\coloneqq\frac{f+1-g}{f+1}\geq\frac1{e}\text{\,?}\]
The elements of $S$ less than $f+1$ are sometimes called sporadic for $S$. Thus $d$ may be considered as the ``density'' of the sporadic elements of $S$. In this note we call $d$ the \emph{Wilf density} of $S$.

The proposed inequality (1) may be considered as a lower bound for the Wilf density in terms of the embedding dimension of $S$. It is well known that (1) holds for $e=2$ \cite{Sylvester} and $e=3$ \cite{Froeberg}.

As an immediate consequence of Zhai's result \cite[Theorem 1]{Zhai}, for $e\geq2$ one obtains the following bounds for the Wilf density.

\noindent\textbf{Proposition}.
\begin{enumerate}
\item $d>\frac16$ for $e=4$ and $d>\frac1{10}$ for $e=5$.
\item If $e\geq2$, then $d\geq\frac2{e^2-e+2}$.
\end{enumerate}
\begin{proof}Let us recall Zhai's result
\noindent\textbf{\cite[Theorem 1]{Zhai}}. If $e\geq2$, then
\[\tag{2}d\geq\frac1e-\frac{m-1}{f+1}\cdot\frac{e-2}{2e}\,.\]
We shall need another

\noindent\textbf{Lemma}. If $e\geq2$, then
\[\tag{3}(e-1)(f+1-g)\geq m-1\,.\]
\noindent\textit{Proof of the Lemma.} Let $L$ be the set of sporadic elements of $S$; hence $\#L=f+1-g$.

Let $m=a_1<a_2<\dots<a_e$ be the atoms of $S$.

If $x\in X\coloneqq\{f+1,f+2,\dots,f+m\}$ is not a multiple of $m$, then there exists an atom $a_i$, $i\geq2$, such that $x-a_i\in L$. Hence
\[\{\,x\in X\mid x\text{ not a multiple of }m\,\}\subseteq Y\coloneqq\bigcup_{i=2, ... ,e}(a_i + L)\,,\]
consequently $m-1\leq\#Y\leq(e-1)(f+1-g)$.\hfill$\qed$Lemma
\vspace{.2cm}\\
\noindent\textit{Continuation of Proof of the Proposition.}
\begin{enumerate}
\item By \cite[Proposition 3.13 and Corollary 6.5]{Eliahou}, in case $f+1\leq3m$, even $d\geq\frac1e$ holds.

In the other case $\frac{f+1}3>m$ ($>m-1$) one can simply replace $m-1$ by $\frac{f+1}3$ in formula (2), this leads to the formula $d>\frac{8-e}{6e}$ (which equals $\frac16$ for $e=4$ and $\frac1{10}$ for $e=5$).
\item From (2) and (3) we get
\[\tag{2*}2(f+1)\leq2e(f+1-g)+(e-2)(m-1)\,.\]
\[\tag{3*}m-1\leq(e-1)(f+1-g)\,.\]
Finally we get b) by substituting (3*) in (2*).
\end{enumerate}
\end{proof}
\end{document}